\documentclass[10pt,a4paper]{article}
\usepackage[utf8]{inputenc}
\usepackage{amsmath}
\usepackage{amsfonts}
\usepackage{amssymb}

\usepackage[utf8]{inputenc}
\usepackage{microtype}

\usepackage{amsthm}

\usepackage{verbatim}
\usepackage{resizegather}
\usepackage{tikz-cd}
\usepackage{babel} 

\usepackage[bottom]{footmisc}

\usepackage{hyperref}

\newtheorem{proposition}{Proposition}[section]
\newtheorem{corollary}{Corollary}[proposition]
\newtheorem{lemma}{Lemma}[proposition]
\newtheorem{theorem}{Theorem}
\newtheorem{remark}{Remark}[proposition]
\newtheorem{example}[proposition]{Example}
\newtheorem{definition}{Definition}

\title{A Formula For the Homology Groups of the Milnor Fibres of an Isolated Real Hypersurface Singularity}
\author{Lars Andersen}
\begin{document}
\maketitle
\begin{abstract}
Given a nonconstant polynomial map $f: \mathbb{R}^{n+1}\to \mathbb{R}$ having an isolated critical point in the origin, and with zero locus of positive dimension, we demonstrate that if any great circle $G_p\subset \mathbb{S}_{\delta}^n$ emancipating at a point $p$ in the sphere Milnor fibre $\mathcal{F}_{\epsilon}^{+}\subset \mathbb{S}_{\delta}^{n}$ has the property that $G_p\cap \left(\mathbb{S}_{\delta}^{n}\setminus \mathcal{F}_{\epsilon}^{+}\right)\neq \emptyset$ then the positive Milnor-Lê fibre $\Phi_{\epsilon}^{+}$ has the homotopy type of the space $\partial \Phi_{\epsilon}^{+}\cup \mathbb{D}^{\lambda(p_1)}\cup\dots\cup\mathbb{D}^{\lambda(p_m)}$, where each disc is attached to the boundary $\partial \Phi_{\epsilon}^{+}$, and where $\lambda(p_i)$ denotes the Morse indices of the critical points of a Morse perturbation of $f$ on the Milnor fibre over the sphere. Under these assumptions a formula for the relative singular homology groups with $\mathbb{Z}$-coefficients $H_{\ast}(\Phi_{\epsilon}^{+}, \partial \Phi_{\epsilon}^{+}; \mathbb{Z})$ is established.
\end{abstract}

\subsection{Introduction}

For isolated complex hypersurface singularities $(\mathbb{C}^{n+1}, 0)\to (\mathbb{C}, 0)$ the Bouquet Theorem of Milnor \cite[Theorem 7.2]{Milnor} implies that the Milnor fibre has homology groups given by
$H_k(\bigvee_{i=1}^{\mu}\mathbb{S}^n; \mathbb{Z})$, for $k\geq 0$,
where $\mu$ is the Milnor number of the singularity. The purpose of this article is to give a similar description of the homology groups of the Milnor fibres in the case of isolated real hypersurface singularities $(\mathbb{R}^{n+1}, 0)\to (\mathbb{R}, 0)$.\\ 

In this case, we will see that under the assumption of a geometric condition which is always valid in low dimensions, a similar formula holds for the homology groups of a real Milnor fibre \emph{relative} to its boundary. The spheres appearing in the bouquet will then no longer be of a fixed dimension, but of variable dimension. The number of spheres, and their dimensions, will be possible to calculate from a certain morsification which we will construct explicitly. This therefore provides an algorithm for computing the homology groups of a Milnor fibre relative to its boundary.\\

To this aim we will \emph{in the large} use ideas which were already present in Milnor's treatise \cite{Milnor} and especially in the proof of the Bouquet Theorem \cite[Theorem 7.2]{Milnor}, adapting Milnor's approach\footnote{Milnor's arguments rely on estimating the Morse indices to be at least $n$ (in fact equality holds as the Milnor fibre has the homotopy type of its complement in the sphere). This estimate does not apply in the case of real singularities.} to the real case by using Stratified Morse Theory as developed by Goresky and MacPherson in \cite{Goresky}, as well as using arguments from Riemannian Geometry.

\subsection{The Local Milnor-Lê Fibration}
Let us be given a non-constant polynomial map germ $f: (\mathbb{R}^{n+1},0)\to (\mathbb{R}, 0)$ such that $(f^{-1}(0),0)$ is of positive dimension and has an isolated singular point in the origin. In this situation there is associated a local Milnor-Lê fibration, and an equivalent fibration over a sphere of sufficiently small radius: 
\begin{theorem}[{\cite[Theorem 13.5]{JS}}]\label{theorem Milnor} There exists a real number $\delta>0$ such that if $0<\epsilon<<\delta$ then 
\begin{equation}\label{ball}
f: f^{-1}([-\epsilon,\epsilon]\setminus\{0\})\cap \bar{\mathbb{B}}_{\delta}^{n+1}\to [-\epsilon,\epsilon]\setminus\{0\}
\end{equation}
is a smooth locally trivial fibration, where $\bar{\mathbb{B}}_{\delta}^{n+1}\subset \mathbb{R}^{n+1}$ is the closed\\
$(n+1)$-dimensional ball centered at the origin and of radius $\delta$. This fibration is equivalent to a fibration
\begin{equation}\label{sphere}
\phi: \mathbb{S}_{\delta}^{n}\setminus \{f=0\}\to \mathbb{S}_{\epsilon}^0
\end{equation}
over the boundary $n$-sphere $\mathbb{S}_{\delta}^n=\partial\bar{\mathbb{B}}_{\delta}^{n+1}$, with projection map $\phi$ such that if $0<\eta<<\epsilon$ and if
$$N(\eta, \delta)=f^{-1}([-\eta, \eta]\setminus \{0\})\cap \mathbb{S}_{\delta}^n$$
then  $\phi_{|N(\eta,\delta)}=f/|f|$.
\end{theorem}

The existence of the first fibration (\hyperref[ball]{\ref*{ball}}) is a standard application of the Thom-Mather First Isotopy Lemma \cite[Proposition 11.1]{Mather}. The existence of the second fibration (\hyperref[sphere]{\ref*{sphere}}) and its equivalence to the first is done by inflating the Milnor tube (see e.g \cite[Lemma 4.4]{Seade}).\\

Let us fix $\delta>0$ and $\epsilon(\delta)>0$ as in Theorem \hyperref[theorem Milnor]{\ref*{theorem Milnor}}. Then
$$K_{\delta}=\{f=0\}\cap \mathbb{S}_{\delta}^n$$
is called the \emph{link} of $f$ at the origin. By transversality of the level hypersurfaces of $f$ and the sphere, $K_{\delta}$ is a smooth $(n-1)$-dimensional manifold. The fibres  
$$\mathcal{F}_{\eta}^{+}=\frac{f}{|f|}^{-1}(\eta)\cap \mathbb{S}_{\delta}^n\setminus K_{\delta},\qquad \mathcal{F}_{\eta}^{-}=\frac{f}{|f|}^{-1}(-\eta)\cap \mathbb{S}_{\delta}^n\setminus K_{\delta}$$
of the second fibration (\hyperref[sphere]{\ref*{sphere}}) are called the \emph{positive} respectively \emph{negative Milnor fibres} of $f$ at the origin. Denoting the fibers of the first fibration (\hyperref[ball]{\ref*{ball}}) by
$$\Phi_{\eta}^{+}=f^{-1}(\eta)\cap \bar{\mathbb{B}}_{\delta}^{n+1},\qquad \Phi_{\eta}^{-}=f^{-1}(-\eta)\cap \bar{\mathbb{B}}_{\delta}^{n+1}$$
there exist by the equivalence of the fibrations (\hyperref[ball]{\ref*{ball}}) and (\hyperref[sphere]{\ref*{sphere}}), stratum-preserving homeomorphisms 
$$\mathcal{F}_{\eta}^{+}\cong \Phi_{\eta}^{+},\qquad \mathcal{F}_{\eta}^{-}\cong \Phi_{\eta}^{-}.$$ 
Since the fibrations are smooth, $\mathcal{F}_{\eta}^{+}$ and $\mathcal{F}_{\eta}^{-}$ are smooth $n$-dimensional manifolds with boundaries.

\subsection{Nondepraved Critical Points} 
Before proceeding we shall need to recall the notion of nondepraved critical points \cite{Goresky} of maps defined on Whitney stratified spaces.\\

Let $\mathbb{K}$ be either the field of real or complex numbers, let $X$ be a $\mathbb{K}$-subanalytic subset of an $\mathbb{K}$-analytic manifold $M$ and let $h: X\to \mathbb{R}$ be the restriction of a smooth function $\tilde{h}$ on the manifold. Let $\mathcal{S}=\{S_i|\ i\in I\}$ be a Whitney stratification of $X$, where the strata are indexed by the partially ordered set $I\subset \mathbb{N}$.\\

In this situation let us recall the following definition.
\begin{definition}[{\cite[Definition 2.1]{Goresky}}] A point $p\in X$ belonging to a stratum $S\in\mathcal{S}$ is called a \emph{critical point} of $h$ if $d\tilde{h}(p)_{|T_p S}=0$. 
\end{definition}

\begin{definition}[{\cite[Definition 2.3]{Goresky}}] The critical point $p\in X$ is \emph{nondepraved} if 
\begin{enumerate} 
\item $p$ is an isolated critical point of $h$.
\item Suppose that $p_i\in S$ converge to $p$, that
$$v_i:=\frac{p_i-p}{|p_i-p|}\to v$$
and that $\ker d h_{|S}(p)\to \tau$. If $v \notin \tau $ then there exists $i_0\in\mathbb{N}\setminus\{0\}$ such that 
$$\forall i\geq i_0,\qquad d h_{|S}(p_i)(v_i)\cdot(h(p_i)-h(p))>0.$$ 
\item For each generalised tangent space $Q$ at $p$, $d\tilde{h}(p)(Q)\neq 0$ except for when $Q=T_p S$.
\end{enumerate}
\end{definition}

In the case $\mathbb{K}=\mathbb{R}$ any critical point of the restriction $h_{|S}$ of $h$ to a stratum $S$ which is furthermore isolated, is nondepraved.

\begin{proposition}[{\cite[Proposition 2.4]{Goresky}}] If $h: \mathbb{R}^n\to\mathbb{R}$ is a real analytic function, and if $p\in\mathbb{R}^n$ is an isolated critical point of $h$, then it is a nondepraved critical point.
\end{proposition}

\subsection{Nonproper Stratified Morse Theory}\label{nonproper}

We shall now recall the main theorems of \emph{nonproper} Stratified Morse Theory. The main reference is \cite[Chapter 10]{Goresky}.\\

Suppose that $M$ is a smooth manifold and let $X\subset M$ be a closed subset endowed with a Whitney stratification $\mathcal{S}=\{S_i| i\in I\}$ with strata indexed by a partially ordered set $I$. Let $\tilde{X}\subset X$ be a subset which is a union of strata.\\

Let $h: X\to\mathbb{R}$ be the restriction of a smooth function on $M$ such that $h$ is proper with a nondepraved critical point $p\in \tilde{X}$. Let $S(p)$ denote the stratum containing $p$.\\

Let us fix a Riemannian metric $(\cdot,\cdot)_M$ on $M$ and if $\rho>0$ let 
$$\bar{\mathbb{B}}_{\rho}(p)=\{x\in M\ |\ (x, p)_M\leq \rho \}$$
denote the closed ball in $M$ with center $p$ and radius $\rho$.\\

By the Whitney conditions $\partial \bar{\mathbb{B}}_{\eta}^M(p)$ intersects any stratum in $\mathcal{S}$ transversally for all sufficiently small $\rho$. One can therefore assume that the restriction of $h$ to $\bar{\mathbb{B}}_{\rho}^M(p)\cap X$ has no critical points except for $p$. Then, as $p$ is a nondepraved critical point, hence is isolated, there exists $\gamma>0$ such that if $s=h(p)$ then $[s-\gamma, s+\gamma]$ contains no critical values except $s$.\\

Fix such $(\rho, \gamma)$ and let us choose a smooth submanifold $N'$ of $M$ such that
\begin{itemize} 
\item $\dim(N')+\dim S(p)=\dim M,$
\item $\forall S'\in\mathcal{S},\qquad S'\pitchfork N',$
\item $N'\cap S=\{p\}.$
\end{itemize}

\begin{remark} It is essential here that the submanifold $N'$ is transverse not just to the strata of $\tilde{X}$, but to $\emph{all}$ strata of $X$.  
\end{remark}

\begin{definition}[{\cite[Definition 10.1, 10.3]{Goresky}}] Let $p\in X$ and let $S\in\mathcal{S}$ be the stratum containing $p$. The \emph{normal slice} through $S$ relative to $\tilde{X}$ at the point $p$ is 
$$N_{\tilde{X}}(p)=\tilde{X}\cap N'\cap \bar{\mathbb{B}}_{\rho}^M(p).$$
\end{definition}

One defines the stratified sets
$$\tilde{X}_{\leq a}=\tilde{X}\cap h^{-1}(-\infty, a],\qquad \tilde{X}_{[a,b]}=\tilde{X}\cap h^{-1}[a, b].$$

\begin{definition}[{\cite[Definition 10.3]{Goresky}}]\label{Morse data}
The Local Morse datum of $h_{|\tilde{X}}$ at the point $p$ is the pair
$$(\mathcal{A}_{\tilde{X}},\mathcal{B}_{\tilde{X}}):=(\bar{\mathbb{B}}_{\eta}^{M}(p)\cap \tilde{X}\cap h^{-1}([s-\gamma, s+\gamma]), \bar{\mathbb{B}}_{\eta}^{M}(p)\cap \tilde{X}\cap h^{-1}(s-\gamma)).$$
The Tangential Morse datum of $h_{|\tilde{X}}$ at $p$ is the pair 
$$(\mathcal{T}_{\tilde{X}},\mathcal{T}'_{\tilde{X}}):=(\bar{\mathbb{B}}_{\eta}^{M}(p)\cap S(p)\cap h^{-1}([s-\gamma, s+\gamma]), \bar{\mathbb{B}}_{\eta}^{M}(p)\cap S(p)\cap h^{-1}(s-\gamma)).$$
The Normal Morse datum of $h_{|\tilde{X}}$ at $p$ is the pair 
$$(\mathcal{N}_{\tilde{X}}, \mathcal{N}'_{\tilde{X}})=(N_{\tilde{X}}(p)\cap h^{-1}([s-\gamma, s+\gamma]), N_{\tilde{X}}(p)\cap h^{-1}(s-\gamma)).$$
\end{definition}

The Local Morse datum of $h_{|\tilde{X}}$ at the nondepraved critical point $p$ specifies the effect on the topology of $\tilde{X}_{\leq a}$ when $a$ passes through the critical value $s$. Namely, the effect is the topological adjunction of $\mathcal{A}_{\tilde{X}}$ along $\mathcal{B}_{\tilde{X}}$.

\begin{theorem}[{\cite[Theorem 10.4]{Goresky}}]\label{Theorem 10.4} For $\gamma>0$ sufficiently small there is a decomposition preserving homeomorphism 
$$\tilde{X}_{\leq \gamma}\cong \tilde{X}_{\leq -\gamma} \cup_{\mathcal{B}_{\tilde{X}}} \mathcal{A}_{\tilde{X}}.$$
\end{theorem}

One can then prove the following Fundamental Theorem of Stratified Morse Theory in the nonproper case, which says that the Local Morse Datum is the product of the Tangential and Normal Morse data.

\begin{theorem}[{\cite[Theorem 10.5]{Goresky}}]\label{Theorem 10.5} There is a decomposition preserving homeomorphism  
$$(\mathcal{A}_{\tilde{X}},\mathcal{B}_{\tilde{X}})\cong (\mathcal{T}_{\tilde{X}},\mathcal{T}'_{\tilde{X}})\times (\mathcal{N}_{\tilde{X}}, \mathcal{N}'_{\tilde{X}})$$ 
$$:=(\mathcal{T}_{\tilde{X}}\times \mathcal{N}_{\tilde{X}}, \mathcal{T}_{\tilde{X}}\times \mathcal{N'}_{\tilde{X}}\cup \mathcal{T}'_{\tilde{X}}\times \mathcal{N}_{\tilde{X}}).$$
\end{theorem}

Define

\begin{definition} The upper and lower half-links of $X$ at the point $p$ with respect to $h_{|\tilde{X}}$ are
$$(l_{\tilde{X}}^{+}(p), \partial l_{\tilde{X}}^{+}(p))=(N_{\tilde{X}}(p)\cap \bar{\mathbb{B}}_{\rho}^M (p)\cap h^{-1}(s+\gamma), N_{\tilde{X}}(p)\cap \partial\bar{\mathbb{B}}_{\rho}^M (p)\cap h^{-1}(s+\gamma)),$$
$$(l_{\tilde{X}}^{-}(p), \partial l_{\tilde{X}}^{-}(p))=(N_{\tilde{X}}(p)\cap \bar{\mathbb{B}}_{\rho}^M (p)\cap h^{-1}(s-\gamma), N_{\tilde{X}}(p)\cap \partial\bar{\mathbb{B}}_{\rho}^M (p)\cap h^{-1}(s-\gamma))$$
\end{definition}

One can prove (see \cite[Proposition 10.7.1]{Goresky}) the following. If the nondepraved critical point $p\in X$ lies in a stratum of $\tilde{X}$ then the normal Morse datum at $p$ has the homotopy type 
$$(\mathcal{N}_{\tilde{X}},\mathcal{N}_{\tilde{X}})\sim (\text{cone}(l_{\tilde{X}}^{-}(p)), l^{-}(p)).$$ 
If $p$ lies in a stratum of $X\setminus\tilde{X}$ (that is, it is a critical point "at infinity") then one can prove (see \cite[Proposition 10.7.2]{Goresky}) that the normal Morse datum at $p$ has the homotopy type 
$$(\mathcal{N}_{\tilde{X}},\mathcal{N}_{\tilde{X}})\sim (l_{\tilde{X}}^{+}(p), \partial l_{\tilde{X}}^{+}(p)).$$

This implies the following theorem

\begin{theorem}[{\cite[Theorem 10.8]{Goresky}}] Suppose that $[a,b]\subset\mathbb{R}$ contains no critical value of the stratified Morse function $h: X\to\mathbb{R}$ except for an isolated critical value $s\in (a,b)$ corresponding to a nondepraved critical point $p$ belonging to a stratum $S$. If $p\in\tilde{X}$ then there is a homotopy equivalence
$$\tilde{X}_{\leq b}\cong \tilde{X}_{\leq a}\cup_{\partial(\mathbb{D}^{\lambda_S(p)}\times\text{cone}(l_{\tilde{X}}^{-}(p)))} (\mathbb{D}^{\lambda_S(p)}\times\text{cone}(l_{\tilde{X}}^{-}(p))).$$
If $p\in X\setminus \tilde{X}$ then there is a homotopy equivalence
$$\tilde{X}_{\leq b}\cong \tilde{X}_{\leq a}\cup_{\partial(\mathbb{D}^{\lambda_S(p)}\times l_{\tilde{X}}^{+}(p))} (\mathbb{D}^{\lambda_S(p)}\times l_{\tilde{X}}^{+}(p))$$
where $\lambda_S(p)$ is the Morse index of $f_{|S}$ at $p$.
\end{theorem}

\begin{remark} The assumption that $h$ is Morse is only needed in order to obtain the homotopy type of the tangential Morse data.
\end{remark}

\subsection{Lemmas} We shall now show that we can perturb $f$ such that the perturbed function has only nondepraved critical points when restricted to the corresponding perturbation of the (positive) Milnor fibre.\\

Let us assume from now onwards that we are given a non-constant polynomial map germ $f: (\mathbb{R}^{n+1},0)\to (\mathbb{R}, 0)$ such that $(f^{-1}(0), 0)$ is of positive dimension and has an isolated singular point in the origin. Take $f: M\to \mathbb{R}$ to be a representative of this germ on an open neighborhood of the origin $M\subset\mathbb{R}^{n+1}$ and take Milnor data $\delta, \epsilon(\delta)>0$ as in the statement of Theorem \hyperref[theorem Milnor]{\ref*{theorem Milnor}}. 

\begin{lemma}\label{main one} There exists $k\in\mathbb{N}\setminus\{0\}$, a neighborhood $U\subset\mathbb{R}^k$ and a polynomial map $F: M \times U \to \mathbb{R}$ such that
\begin{enumerate}
\item if $f_t=F(x,\cdot)$ then $f_0=f$,
\item there exists an open dense subset $V\subset U$ such that if $t\in V$ then
$$f_{t|\mathbb{S}_{\delta}^n}: \mathbb{S}_{\delta}^n\to \mathbb{R}$$ has only nondegenerate critical points with pairwise distinct critical values. Moreover the critical values of $f_t$ belong to $(-\epsilon,\epsilon)$.
\end{enumerate}
\end{lemma}
\begin{proof} Consider the polynomial map
$$F: M\times \mathbb{R}^{n+1}\to \mathbb{R},\qquad F(x, t)=f(x)-\sum_{i=1}^{n+1} t_i x_i.$$
where $t=(t_1,\dots, t_{n+1})$ are coordinates on $\mathbb{R}^{n+1}$ and where $x=(x_1,\dots, x_{n+1})$ are coordinates on the neighborhood $M$. Let 
$$\Psi: M\times \mathbb{R}^{n+1}\to \mathbb{R}^{n+1},$$
$$\Psi(x,t)=df_t(x).$$
Then the Jacobian of $\Psi$ is given by the $(n+1)\times (2n+2)$-matrix
\[
\text{Jac}(\Psi)(x,t)
=
\begin{bmatrix}
    \partial^2 f/\partial x_1^2 & \dots & \partial^2 f/\partial x_1\partial x_{n+1} & 1  & \dots &  0\\
    \vdots & \dots & \vdots & \vdots & \ddots & \vdots \\
    \partial^2 f/\partial x_1\partial x_{n+1} & \dots & \partial^2 f/\partial x_{n+1}^2 & 0 & \dots  & 1
\end{bmatrix}
\]
This matrix is of maximal rank $n+1$ so $\Psi$ is a submersion. One can therefore apply \cite[Theorem 2.2.3]{Goresky} to $F$ and to the compact subset $\mathbb{S}_{\delta}\subset \mathbb{R}^{n+1}$ to deduce that there exists a dense open subset $\tilde{V}\subset \mathbb{R}^{n+1}$ of the parameter space such that if $t\in \tilde{V}$ then
$$f_{t| \mathbb{S}_{\delta}}: \mathbb{S}_{\delta}\to \mathbb{R}$$ 
has only nondegenerate critical points with pairwise disctinct critical values. Since the origin is an isolated critical point of $f$, the critical points of $f_t: M\to\mathbb{R}$ tends to the origin as $|t|\to 0$ so there exists $\tilde{t}_0\in\mathbb{R}$ such that the critical values of $f_t$ belong to $(-\epsilon,\epsilon)$ whenever $t\in \mathbb{R}^{n+1}$ is such that $|t|< \tilde{t}_0$. One takes   
$$U=\mathbb{R}^{n+1}\cap \{|t|< \tilde{t}_0\},\qquad V=\tilde{V}\cap U$$ 
and then the restricted map $F: M\times U\to \mathbb{R}$ satisfies the wanted properties.
\end{proof}

Let us choose $k, F, U,$ and $V$ as in the previous Lemma \hyperref[main one]{\ref*{main one}}. The proof of the following lemma is inspired by the proof of \cite[Lemma 4.4]{Seade}.

\begin{lemma}\label{main two} There exists $t_0\in \mathbb{R}$ such that for any $t\in V$ such that $|t|\leq t_0$ there exists a stratum-preserving homeomorphism 
$$\psi_t: \mathcal{F}_{\epsilon}^{+}\to \phi_t^{-1}(\epsilon)\cap\{f_t\geq\epsilon\}\cap \mathbb{S}_{\delta}^n$$
where 
$\phi_t: \mathbb{S}_{\delta}\cap \{f_t\geq \epsilon\}\to \mathbb{R}$ 
is a smooth function such that $\phi_0=\phi$ where $\phi$ is the map in Theorem \hyperref[theorem Milnor]{\ref*{theorem Milnor}}.
\end{lemma}
\begin{proof} By Lemma \hyperref[main one]{\ref*{main one}} the critical values of $f_t$ tends to the origin as $|t|\to 0$. Therefore there exists $t_0\in\mathbb{R}$ such that for any $t\in V$ of modulus $|t|\leq t_0$ the critical values of $f_t: M\to \mathbb{R}$ all lie in the interval $(-\epsilon, \epsilon)$. Fix $t\in V$ such that $|t|\leq t_0$ and let $\eta>0$ be a real number such that $[\eta, \epsilon)$ contains no critical values of $f_t$.\\  

Let $\textbf{x}$ be the restriction of the normal vector field on $\bar{\mathbb{B}}_{\delta}$ to $\bar{\mathbb{B}}_{\delta}\cap\{f_t\geq\eta\}$. By the choice of the Milnor radius $\delta$ it follows that the integral lines of the normal vector field are transverse to each sphere contained in $\bar{\mathbb{B}}_{\delta}$ hence this is also true for the integral lines of the restriction $\textbf{x}$. Moreover $f_t^{-1}(s)$ is transverse to the integral lines of the gradient $\text{grad}(f_t)$ as long as $s\geq \eta$.\\

Let us define a vector field $\textbf{v}_t$ on $\bar{\mathbb{B}}_{\delta}\cap\{f_t\geq \eta\}$ by
$$\textbf{v}_t=\lVert \textbf{x} \rVert \text{grad}(f_t)+\lVert \text{grad}(f_t) \rVert \textbf{x}.$$
Then $\textbf{v}_t$ bisects the angle between $\mathbf{x}$ and $\text{grad}(f_t)$. By Cauchy-Schwartz
$$\textbf{v}_t\cdot\mathbf{x}=\lVert \textbf{x} \rVert\left(\text{grad}(f_t)\cdot \mathbf{x}+\lVert\text{grad}(f_t)\rVert \lVert\mathbf{x}\rVert\right)$$
$$\geq \lVert \textbf{x} \rVert\left(\text{grad}(f_t)\cdot \mathbf{x}+\lVert\text{grad}(f_t)\cdot \mathbf{x}\rVert\right)\geq 0$$
and similarly 
$$\textbf{v}_t\cdot\text{grad}(f_t)\geq 0.$$
The vector fields $\text{grad}(f_t)$ and $\mathbf{x}$ can be identified with the gradient vector fields of non-negative functions, namely 
$$g(x)=f^2(x),\quad \text{and}\quad h=\sum_{i=1}^{n+1} x_i^2$$ 
respectively. Since $\textbf{v}_t$ bisects these vector fields it follows from \cite[Corollary 3.4]{Milnor} that $\textbf{v}_t\cdot\mathbf{x}=0$ at a point $p\in f_t^{-1}(s)$ if and only if either $\mathbf{x}$ is zero at $p$, or $\text{grad}(f_t)(p)=0$.\\

In the first case we get a contradiction because if $t_0$ is sufficiently small then the origin is not contained in any fibre $f_t^{-1}(s)$ with $s\geq\eta$. In the second we get a contradiction because $f_t^{-1}(s)$ is nonsingular. Therefore
\begin{equation}\label{inc}
\textbf{v}_t\cdot \textbf{x}>0,\qquad \textbf{v}_t\cdot \text{grad}(f_t)>0.
\end{equation}
with strict inequalities. Since $\bar{\mathbb{B}}_{\delta}\cap\{f_t\geq \eta\}$ is compact the vector field $\textbf{v}_t$ can be integrated, by the Picard-Lindelöf theorem. One defines a stratum-preserving homeomorphism 
$$h_t: f_t^{-1}([\eta, \epsilon])\cap\bar{\mathbb{B}}_{\delta}\to \mathbb{S}_{\delta}\cap \{f_t\geq \eta\}$$
by sending $x$ to the point $h_t(x)$ where the integral line of $\textbf{v}_t$ starting at $x$ intersects the sphere $\mathbb{S}_{\delta}$. This point exists because the gradient and the normal vectors strictly increases along the integral lines of $\textbf{v}_t$ by (\hyperref[inc]{\ref*{inc}}).\\

One then defines
$$\phi_t: \mathbb{S}_{\delta}\setminus \{f_t\geq \eta\}\to \mathbb{R},\qquad \phi_t=f_t\circ h_t^{-1}.$$
In particular
\begin{equation}\label{ettan}
\phi_t^{-1}(\epsilon)\cap \mathbb{S}_{\delta}\cap \{f_t\geq \epsilon\}\cong f_t^{-1}(\epsilon)\cap\bar{\mathbb{B}}_{\delta}.
\end{equation}
Now since $\epsilon$ is not a critical value of $f_{t| \mathbb{B}_{\delta}}$ for any $t\in V$ such that $|t|\leq t_0$, and since $f_t^{-1}(\epsilon)$ is transverse to the sphere $\mathbb{S}_{\delta}$ a standard argument using the First Isotopy Lemma of Thom-Mather \cite[Proposition 11.1]{Mather} shows that there exists a stratum-preserving homeomorphism
\begin{equation}\label{tvaan}
f_t^{-1}(\epsilon)\cap\bar{\mathbb{B}}_{\delta}\cong f^{-1}(\epsilon)\cap\bar{\mathbb{B}}_{\delta}.
\end{equation}
Composing the homeomorphisms (\hyperref[ettan]{\ref*{ettan}}) and (\hyperref[tvaan]{\ref*{tvaan}}) gives the wanted homeomorphism $\psi_t$.
\end{proof}

In what follows we shall for brevity write $V$ instead of $V\cap\{|t|\leq t_0\}$. Furthermore let us define
\begin{definition} For any $t\in V$ let
$$\mathcal{F}_t^{+}:=\phi_t^{-1}(\epsilon)\cap\{f_t\geq \epsilon\}\cap\mathbb{S}_{\delta}^n$$
\end{definition}

Fix $t\in V$ from now onwards. Stratify $\mathcal{F}_t^{+}$ with strata 
$$S=\phi_t^{-1}(\epsilon)\cap\{f_t>\epsilon\}\cap \mathbb{S}_{\delta}^n$$
$$S'=\phi_t^{-1}(\epsilon)\cap f_t^{-1}(\epsilon)\cap\mathbb{S}_{\delta}^n$$
and let us remark that by the previous lemma, $S$ is a smooth manifold of dimension $n$, homeomorphic to
$$f^{-1}(\epsilon)\cap\mathbb{B}_{\delta}$$
and that $S'$ is a smooth manifold of dimension $n-1$, homeomorphic to the boundary
$$f^{-1}(\epsilon)\cap \mathbb{S}_{\delta}$$
The idea is now to consider the restriction of $f_t$ to $\mathcal{F}_t^{+}$ and apply Stratified Morse Theory to this function. To begin with we can analyse the critical points in the interior stratum.

\begin{lemma} Let $t\in V.$ The set of critical points of the restriction of $f_t$ to the stratum $S$ is
$$\mathcal{C}=\{x\in S\ |\ \exists \lambda\in\mathbb{R}\quad \text{grad}(f_t)(x)=\lambda x\}$$
 
\end{lemma}
 \begin{proof} A point $x\in S$ is a critical point of $f_{t|S}$ if and only if the vector $x$ is orthogonal to the gradient,
 $$(x, \text{grad}(f_{t})(x))=0$$
Since $S=\text{int}\mathcal{F}_t^{+}$ is an open subset of the $n$-sphere $\mathbb{S}_{\delta}$ of dimension $n$, we can identify the set of normal vectors to $S$ at a point $x$ with the set of normal vectors to $\mathbb{S}_{\delta}$ at $x$. The former is spanned by the vector $x$, and therefore  $(x, \text{grad} (f_{t})(x))=0$ if and only if $\text{grad}(f_{t})(x)$ is a real multiple of $x$.
 \end{proof}
 
\begin{lemma}\label{main three} The critical points of $f_{t|S}$ are nondegenerate.
\end{lemma} 
\begin{proof} Since the interior of $\mathcal{F}_t^{+}$ is an open subset of the $n$-sphere $\mathbb{S}_{\delta}$ of dimension $n$, and since being a nondegenerate critical point is independent of the choice of coordinates (see e.g \cite{Morse}), this follows from the fact that $f_{t| \mathbb{S}_{\delta}}$ has only nondegenerate critical points.
\end{proof}
\begin{remark} In particular the critical points of $f_{t|S}$ are isolated hence nondepraved.
\end{remark} 

\subsection{The Main Result} 
We retain the notation of the previous subsection. Let $\delta>0$ and $\epsilon(\delta)>0$ be as in Theorem \hyperref[theorem Milnor]{\ref*{theorem Milnor}}. Let $k(\epsilon)\in \mathbb{N}\setminus\{0\}$, the neighborhood $U(\epsilon)\subset \mathbb{R}^k$, the subset $V(\epsilon)\subset U(\epsilon)$ and the polynomial map $F(\epsilon)$ be as given by Lemma \hyperref[main one]{\ref*{main one}}. Let us furthermore assume that $V(\epsilon)$ is chosen such that the conclusion of \hyperref[main two]{\ref*{main two}} is valid whenever $t\in V(\epsilon)$. Fix each one of these entities (including $t\in V$). Recall that 
$$\mathcal{F}_{t, \epsilon}^{+}=\phi_t^{-1}(\epsilon)\cap \{f_t\geq \epsilon\}\cap\mathbb{S}_{\delta}^n$$ 
denotes the perturbation of the Milnor fibre corresponding to $f_t=F(\cdot, t)$.

\begin{theorem}\label{goal} Let $f: (\mathbb{R}^{n+1}, 0)\to (\mathbb{R}, 0)$ be a nonconstant polynomial map germ such that $(f^{-1}(0), 0)$ is of positive dimension and has an isolated singular point in the origin. Suppose that for any $p\in \mathcal{F}_{\epsilon}^{+}$ if $G_p\subset \mathbb{S}_{\delta}$ is a great circle containing $p$ then $G_p\cap \mathbb{S}_{\delta}\setminus \mathcal{F}_{\epsilon}^{+}\neq\emptyset$. Then there exists a homotopy equivalence
$$\Phi_{\epsilon}^{+}\sim \partial \Phi_{\epsilon}^{+}\cup  (\bigcup_{i=1}^m \mathbb{D}^{\lambda(p_i)}),$$ 
where $m$ is the number of critical points $p_i$ of $f_{t| \text{int}\mathcal{F}_{t, \epsilon}^{+}}$, where $\lambda(p_i)$ denotes the Morse indices, and where each disc is attached via attaching maps 
$$r_i: \partial\mathbb{D}^{\lambda(p_i)}\to \partial \Phi_{\epsilon}^{+},\qquad i=1,\dots, m.$$
\end{theorem}
\begin{proof} We will apply Stratified Morse Theory \cite{Goresky} to $M=\mathbb{S}_{\delta}^{n}\subset \mathbb{R}^{n+1}$, the closed subset $X=\mathcal{F}_{t, \epsilon}^{+}\subset \mathbb{S}_{\delta}$ with the Whitney stratification $\mathcal{S}=\{S, S'\}$, the subset $\tilde{X}=S$ and the proper smooth function 
$$f_{t| \mathcal{F}_{t, \epsilon}^{+}}: \mathcal{F}_{t, \epsilon}^{+}\to \mathbb{R}.$$

Let $p_1,\dots, p_m\in S$ denote the critical points of $f_{t| S}$ and let $s_i=f_t(p_i)$ be their values, arranged so that $s_1<s_2<\dots<s_m$. By Lemma \hyperref[main three]{\ref*{main three}} the critical points $p_i$ are nondegenerate;  let $\lambda(p_i), i=1,\dots, m$ denote the Morse indices.\\

By compactness of the sphere $\mathbb{S}_{\delta}$ it follows from the Curve Selection Lemma (see e.g \cite[Lemma 5.7]{Milnor}) that the origin in $\mathbb{R}$ is not an accumulation point of critical values of $f_{t| S}$. We can therefore assume that $t\in V$ is chosen sufficiently small in order to have that $s_i>\epsilon$ for $i=1,\dots, m$.\\
 
For each $i=1,\dots, m$ choose data $(\eta_i, \gamma_i)$ in accordance with the procedure in section \hyperref[nonproper]{\ref*{nonproper}} and let
$$(\mathcal{N}_{S}(p_i), \mathcal{N}'_{S}(p_i))=(N_{S}(p_i)\cap f_t^{-1}([s_i-\gamma_i, s_i+\gamma_i]), N_{S}(p_i)\cap f_t^{-1}(s_i-\gamma_i))$$
denote the Normal Morse datum (Definition \hyperref[Morse data]{\ref*{Morse data}}) of $f_{t|S}$ at $p_i$. Since $S$ is a stratum of maximal dimension $\dim S=\dim \mathbb{S}_{\delta}$ it follows that for any $i=1,\dots, m$ the normal slice is 
$$N_{S}(p_i)=S\cap N'(p_i)\cap \bar{\mathbb{B}}_{\eta_i}^M (p_i)=\{p_i\}.$$ Indeed the corresponding transverse submanifold $N'(p_i)$ is zero dimensional and contains $p_i$ so its intersection with $\bar{\mathbb{B}}_{\eta_i}^{M}(p_i)$ is reduced to $p_i$, if $\eta_i$ is sufficiently small.\\

It follows that the Normal Morse datum is
$$(\mathcal{N}_{S}(p_i), \mathcal{N}'_{S}(p_i))=(p_i, \emptyset),\qquad i=1,\dots, m$$ 
and so by the Main Theorem \hyperref[Theorem 10.5]{\ref*{Theorem 10.5}} the Local Morse datum of $f_{t| S}$ at $p_i$ is homeomorphic by a stratum-preserving homeomorphism to the Tangential Morse datum 
$$(\mathcal{T}_{S}(p_i),\mathcal{T}'_{S}(p_i))=(\bar{\mathbb{B}}_{\eta_i}^{M}(p_i)\cap S(p_i)\cap f_t^{-1}([s_i-\gamma_i, s+\gamma_i]), \bar{\mathbb{B}}_{\eta_i}^{M}(p_i)\cap S(p_i)\cap f_t^{-1}(s_i-\gamma_i))$$
of $f_{t|S}$ at $p_i$. Since the point $p_i$ is nondegenerate by Lemma \hyperref[main three]{\ref*{main three}} and since $\dim S=n$, there exists by \cite[Proposition 4.5]{Goresky} a homeomorphism of pairs 
\begin{equation}
(\mathcal{T}_{S}(p_i),\mathcal{T}'_{S}(p_i)\to (\mathbb{D}^{\lambda(p_i)}\times \mathbb{D}^{n-\lambda(p_i)}, \partial \mathbb{D}^{\lambda(p_i)}\times \mathbb{D}^{n-\lambda(p_i)})
\end{equation}
By Theorem \hyperref[Theorem 10.5]{\ref*{Theorem 10.5}} \cite[Theorem 10.5]{Goresky} one concludes that there exists a stratum-preserving homeomorphism
\begin{equation}\label{eq one}
S\cong S\cap f_t^{-1}([\epsilon, s_1-\gamma_1])\cup \left(\bigcup_{i=1}^m \mathbb{D}^{\lambda(p_i)}\times \mathbb{D}^{n-\lambda(p_i)} \right) 
\end{equation}
where each $\lambda(p_i)$-handle is attached via attaching maps 
\begin{equation}\label{darn}
\partial \mathbb{D}^{\lambda(p_i)}\times \mathbb{D}^{n-\lambda(p_i)}\to S\cap f_t^{-1}([\epsilon, s_i-\gamma_i]).
\end{equation}
To complete the proof we shall proceed in several steps.
\begin{enumerate} 
\item Each $\lambda(p_i)$-handle is attached to $S\cap f_t^{-1}(s_i-\gamma_i)$. Indeed each of the attaching maps in (\hyperref[darn]{\ref*{darn}}) is by the proof of \cite[Proposition 4.5]{Goresky} given by the composition between a stratum-preserving homeomorphism 
$$\phi_i': \partial \mathbb{D}^{\lambda(p_i)}\times \mathbb{D}^{n-\lambda(p_i)}\to \bar{\mathbb{B}}_{\eta_i}^{M}(p_i)\cap S\cap f_t^{-1}(s_i-\gamma_i)$$
and the inclusion 
$$i: \bar{\mathbb{B}}_{\eta_i}^{M}(p_i)\cap S\cap f_t^{-1}(s_i-\gamma_i)\hookrightarrow S\cap f_t^{-1}(s_i-\gamma_i)$$
so the handles are attached as claimed.
\item We claim that there exists for each $i=1,\dots, m$ a continuous map
$$\phi_i: \bar{\mathbb{B}}_{\eta_i}^{M}(p_i)\cap S\cap f_t^{-1}(s_i-\gamma_i)\to f_t^{-1}(\epsilon)\cap \mathbb{S}_{\delta}$$
which is a homeomorphism onto its image. Note that by construction 
$$L_i:=\bar{\mathbb{B}}_{\eta_i}^{M}(p_i)\cap S\cap f_t^{-1}(s_i-\gamma_i)=\bar{\mathbb{B}}_{\eta_i}^i(p_i)\cap \mathbb{S}_{\delta}\cap f_t^{-1}(s_i-\gamma_i).$$
For each $i=1,\dots, m$ consider the normal vector field $\mathbf{x}_{p_i}$ on $\mathbb{S}_{\delta}$ centered at the point $p_i$. For each $x\in L_i$ choose a tangent vector $\vec{v}_x\in T_x\mathbb{S}_{\delta}$ parallel to the flow of $\mathbf{x}_{p_i}$ at $x$ and let $\gamma_{x}(t)$ be the geodesic flow in the direction $\dot{\gamma}_{x}(0)=\vec{v}_x$ starting at $\gamma_{x}(0)=x$. Define 
$$\phi_i: L_i\to \mathbb{S}_{\delta}\cap f_t^{-1}(\epsilon)$$
by sending $x\in L_i$ to the point $\phi_i(x)=\gamma_{x}(t_0)$ where the geodesic flow first intersects $\mathbb{S}_{\delta}\cap f_t^{-1}(\epsilon)$. Indeed the point $\phi_i(x)$ exists because
$$\gamma_{x}: \mathbb{R}\to \mathbb{S}_{\delta},\qquad t\mapsto \gamma_{x}(t)$$
defines a great circle emancipating at $p_i$ and by hypothesis this great circle is not entirely contained in $\mathcal{F}_{\epsilon}^{+}$ hence not in $\mathcal{F}_{t,\epsilon}^{+}$ if $t_0$ is chosen sufficiently small.  
\begin{enumerate} 
\item $\phi_i$ is continuous with continuous inverse. Indeed by the Picard-Lindelöf Theorem the choice of a geodesic on a Riemannian Manifold through a given point $x$ in a direction $\vec{v}$ varies differentiably with respect to $x$ and $\vec{v}$.
\item $\phi_i$ is injective. Indeed suppose $x\neq x'$ but $\phi_i(x)=\phi_i(x')$. This means that there exists $t_0, t_0'\in\mathbb{R}$ such that
$$\gamma_{x}(t_0)=\gamma_{x'}(t_0')$$
so that the geodesics $\gamma_x: t\mapsto \gamma_x(t)$ and $\gamma_{x'}: t\mapsto \gamma_{x'}(t)$ intersect. Since $\vec{v}_x$ and $\vec{v}_{x'}$ are parallel to the flow of the normal vector field in $\mathbb{S}_{\delta}$ centered at $p_i$ they are parallel to great spheres emancipating at $p_i$. Since $\gamma_x$ and $\gamma_x$ are geodesics they parallel transport $\vec{v}_x$ and $\vec{v}_{x'}$, respectively. Therefore $\vec{v}_x$ and $\vec{v}_{x'}$ are either parallel or $x,x'\in L_i$ are antipodal points on $\mathbb{S}_{\delta}$, which is impossible if $\eta_i$ is chosen small enough. Therefore the vectors $\vec{v}_x$ and $\vec{v}_{x'}$ are parallel. By the Morse Lemma \cite[Lemma 2.2]{Morse} one can choose local coordinates on $\mathbb{S}_{\delta}\cap \bar{\mathbb{B}}_{\eta_i}^{M}(p_i)$ in which $p_i$ is the origin and in which
$$f_{t| \mathbb{S}_{\delta}\cap \bar{\mathbb{B}}^{M}_{\eta_i}(p_i)}=\sum_{i=1}^{\lambda(p_i-1)} x_i^2-\sum_{j=\lambda(p_i)}^{n} x_j^2$$
That $x,x'\in L_i$ and that $\vec{v}_{x}$ and $\vec{v}_{x'}$ are parallel means with this choice of coordinates that $x,x'$ lie on the same ray $l=\{\tau\vec{w}\ |\ \tau\in \mathbb{R}^{>}\}$ through the origin. Write 
$$x=(q_1,\dots, q_n),\qquad \vec{w}=\vec{x}.$$
Then $x,x'\in l\cap L_i$ if and only if
\begin{equation}\label{barn}
\sum_{i=1}^{\lambda(p_i-1)} (\tau q_i)^2-\sum_{j=\lambda(p_i)}^{n} (\tau q_j^2)=\gamma_i.
\end{equation}
But (\hyperref[barn]{\ref*{barn}}) has at most one solution for $\tau>0$. Therefore $x=x'$. This contradiction shows that $\phi_i$ is injective. 

\end{enumerate}
It follows that $\phi_i: L_i\to \phi_i(L_i)$ is a homeomorphism and one defines 
$$\tilde{r}_i=h_i\circ \phi_i\circ\phi_i': \partial \mathbb{D}^{\lambda(p_i)}\times \mathbb{D}^{n-\lambda(p_i)}\to \partial \Phi_{\epsilon}^{+}$$
where $h_i$ is as in the proof of Lemma \hyperref[main two]{\ref*{main two}}.
\item We claim that there exists a homotopy equivalence
$$S\cap f_t^{-1}([\epsilon, s_1-\gamma_1])\sim \partial \Phi_{\epsilon}^{+}$$
To establish this one considers the function 
$$\alpha_t=f_t-\epsilon: S\cap f_t^{-1}([\epsilon, s_1-\gamma_1])\to \mathbb{R}.$$
To simplify the notation, write $s=s_1-\gamma_1$. The following properties are immediately verified:
\begin{enumerate}
\item $S\cap f_t^{-1}([\epsilon, s])$ is compact and semialgebraic,
\item $\alpha_t$ is proper and semialgebraic, 
\item $\alpha_t(x)\geq 0$ for all $x\in S\cap f_t^{-1}([\epsilon, s])$ and $\alpha_t^{-1}(0)=S\cap f_t^{-1}(\epsilon)$.
\end{enumerate}

In the notation of \cite{durfee}, $\alpha_t$ is a therefore a \emph{semialgebraic rug function} relative to $S\cap f_t^{-1}(\epsilon)\subset S\cap f_t^{-1}([\epsilon, s])$. Write $s=\epsilon+\rho$ for some $\rho>0$. Then 
$$\alpha_t^{-1}([0, \rho])=S\cap f_t^{-1}([\epsilon, s])$$
is a semialgebraic neighborhood of $S\cap f_t^{-1}(\epsilon)$. In particular if $\rho$ is small enough, then by \cite[Proposition 1.6]{durfee} the inclusion 
$$f_t^{-1}(\epsilon)\cap S\to S\cap f_t^{-1}([\epsilon, s])$$
is a homotopy equivalence, which establishes the claim. 
\end{enumerate}

The inclusion 
$$j: (\mathbb{D}^{\lambda}, \partial \mathbb{D}^{\lambda})\to (\mathbb{D}^{\lambda}\times \mathbb{D}^{n-\lambda}, \partial \mathbb{D}^{\lambda}\times \mathbb{D}^{n-\lambda})$$ 
is a homotopy equivalence so, using \cite[Lemma 3.7]{Morse}, one obtains a homotopy equivalence
$$S\sim \partial \Phi_{\epsilon}^{+}\cup_{\partial \mathbb{D}^{\lambda(p_1)}} \mathbb{D}^{\lambda(p_1)}\cup\dots\cup_{\partial \mathbb{D}^{\lambda(p_m)}} \mathbb{D}^{\lambda(p_m)}$$
with attaching maps $r_i:= j\circ \tilde{r}_i$. But this is what was to be proven because the manifold with boundary $\Phi_{\epsilon}^{+}$ has by the Tubular Neighborhood Theorem \cite{Hatcher} the homotopy type of its interior which is by Lemma \hyperref[main two]{\ref*{main two}} homeomorphic to $S$.

\end{proof}
\begin{remark} The hypothesis in Theorem \hyperref[goal]{\ref*{goal}} that any great circle $G_p$ emancipating at a point $p$ in the interioir of the Milnor fibre is not entirely contained in it, is always fulfilled if $n=1$ because then the boundary of the Milnor fibre $\mathcal{F}_{\epsilon}^{+}\subset \mathbb{S}_{\delta}^1$ is a smooth zero-dimensional submanifold of the circle hence a disjoint union of points bounding the arcs $\Gamma_1,\dots, \Gamma_M$ constituting the interior of the Milnor fibre. In the case $n=2$ the hypothesis is likewise fulfilled because then the boundary of $\mathcal{F}_{\epsilon}^{+}\subset \mathbb{S}_{\delta}^2$ is a one-dimensional smooth manifold hence a disjoint union of closed curves without self-intersections, by the classification of compact one-dimensional smooth manifolds. Each curve separates the sphere $\mathbb{S}_{\delta}^2$ into two bounded parts, hence the claim.
\end{remark}
\subsection{Concluding Remarks}

As a corollary one immediately deduces the following result.

\begin{corollary}\label{corr} Keep the notation of Theorem \hyperref[goal]{\ref*{goal}} and suppose that $\lambda(p_i)>0$ for each $i=1,\dots, m$. The singular homology groups of the positive Milnor fibre of $f$ relative to its boundary satisfy
$$H_k(\Phi_{\epsilon}^{+}, \partial \Phi_{\epsilon}^{+}; \mathbb{Z})=\tilde{H}_{k-1}(\bigvee_{i=1}^m \mathbb{S}^{\lambda(p_i)-1}; \mathbb{Z}),\qquad k\geq 1.$$
\end{corollary}
\begin{proof} By Theorem \hyperref[goal]{\ref*{goal}}, there is a homotopy equivalence
$$S\sim \partial \Phi_{\epsilon}^{+}\cup_{\partial \mathbb{D}^{\lambda(p_1)}} \mathbb{D}^{\lambda(p_1)}\cup\dots\cup_{\partial \mathbb{D}^{\lambda(p_m)}} \mathbb{D}^{\lambda(p_m)}$$
where each disc $\mathbb{D}^{\lambda(p_i)}$ attached along its boundary sphere via attaching maps
$$j\circ h_t^{-1}\circ i\circ H_i: \partial \mathbb{D}^{\lambda(p_i)}\to \partial \Phi_{\epsilon}^{+}.$$
Applying the Excision Theorem \cite{Hatcher} to the complement of the images of these maps one obtains

$$H_k(\Phi_{\epsilon}^{+}, \partial \Phi_{\epsilon}^{+}; \mathbb{Z})=\bigoplus_{i=1}^m H_k(\mathbb{D}^{\lambda(p_i)}, \partial \mathbb{D}^{\lambda(p_i)}; \mathbb{Z})$$ 
Using the long exact sequences of the pairs $(\mathbb{D}^{\lambda(p_i)}, \partial \mathbb{D}^{\lambda(p_i)})$ gives
$$H_k(\mathbb{D}^{\lambda(p_i)}, \partial \mathbb{D}^{\lambda(p_i)}; \mathbb{Z})=\tilde{H}_{k-1}(\mathbb{S}^{\lambda(p_i)-1}; \mathbb{Z})$$
hence 
$$H_k(\Phi_{\epsilon}^{+}, \partial \Phi_{\epsilon}^{+}; \mathbb{Z})=\bigoplus_{i=1}^m \tilde{H}_{k-1}(\mathbb{S}^{\lambda(p_i)-1}; \mathbb{Z})$$
$$=\tilde{H}_{k-1}(\bigvee_{i=1}^m \mathbb{S}^{\lambda(p_i)-1}; \mathbb{Z})$$ 
\end{proof}

We now give a example, illustrating the fact that finding the Morse indices in our case is a less trivial matter than in the case of complex analytic singularities.

\begin{example} Consider the cusp singularity $f: \mathbb{R}^2\to \mathbb{R}$ given by $f(x,y)=x^3-y^2$. Let $\epsilon>0$ be such that 
$$f: f^{-1}((0, \epsilon])\cap \bar{\mathbb{B}}^2\to (0, \epsilon]$$
is a trivial fibration, where $\mathbb{B}^2=\mathbb{B}_1^2$ is the unit ball centered at the origin in $\mathbb{R}^2$. Let
$$\Phi_{\eta}^{+}=f^{-1}(\epsilon)\cap \mathbb{B}^2$$ 
$$\mathcal{F}_{\epsilon}^{+}=\{\frac{f}{|f|}=\epsilon,\quad f\neq 0,\quad x^2+y^2=1 \}$$
be positive Milnor fibres over the ball and over the sphere, respectively. Take 
$$F: \mathbb{R}^2\times \mathbb{R}\to \mathbb{R},\qquad F(x,t)=f+3tx$$
and define $\tilde{f}=f_{1| \mathbb{S}^1}$. Take local coordinates on $\mathcal{F}_{\epsilon}^{+}$ given by stereographical projection:
$$\phi: \Phi_{\epsilon}^{+}\to \mathbb{R},\qquad \phi(x,y)=\frac{x}{1-y}$$
Then 
$$\phi^{-1}(z)=(\frac{2z}{z^2+1}, \frac{z^2-1}{z^2+1})$$
so that the "height function" $\tilde{f}$ is expressed in this chart as 
$$\tilde{f}\circ \phi^{-1}(z)=\frac{1}{z^2+1}\left(6z+\frac{8z^3}{(z^2+1)^2}-\frac{(z^2-1)^2}{z^2+1}\right)$$
and one finds that there is one critical point $p=(1, 0)\in \mathcal{F}_{\epsilon}^{+}$ corresponding to $z=1$. The second derivative of $(\tilde{f}\circ \phi)^{-1}$ is 
$$(\tilde{f}\circ \phi^{-1})''(z)=u(z)(-8+60 z-120 z^3+40 z^4-36 z^5)$$
where $u: \mathbb{R}\to \mathbb{R}$ a positive function. Evaluating at $\phi^{-1}(p)=1$ one finds 
$$(\tilde{f}\circ \phi^{-1})''(1)=-24.$$ 
Thus $\tilde{f}$ is Morse at $p$ and the Morse index is $\lambda(p)=1$. According to Theorem \hyperref[goal]{\ref*{goal}} it follows that 
$$\Phi_{\epsilon}^{+}\sim \partial \Phi_{\epsilon}^{+}\cup \mathbb{D}^1.$$
In this case the boundary $\partial \Phi_{\epsilon}^{+}$ consists of two disjoint points. It follows that the Milnor fibre has the homotopy type of a line segment, hence is contractible. 
\end{example} 

In the complex case (see e.g \cite{Ebeling}) the intersection matrix of a set of strongly distinguished vanishing cycles determines the monodromy action on the Milnor fibre, and the intersection pairing satisfies the Picard-Lefschetz formula. In the real case, the monodromy is trivial since the fibration 
$$f: f^{-1}((0,\epsilon])\cap \bar{\mathbb{B}}_{\delta}\to (0, \epsilon]$$
is trivial. It would therefore be a very interesting problem to find the intersection matrix of the generators of the homology group in Corollary \hyperref[corr]{\ref*{corr}}.

\bibliographystyle{plain}
\bibliography{biblio}

\end{document}